%% file: ncube3fin.tex
\def\version{May 6, 2003}

\documentclass[12pt]{article}
\usepackage[psamsfonts]{amsfonts}
\usepackage{amsmath,amssymb,amsthm}
\usepackage[dvips]{graphicx}
\usepackage{eepic}

\input def99c 
\UseSection  
\renewcommand{\to}      {\rightarrow}

\newcommand{\Pro}{{\mathbb P}_p}
\newcommand{\Exp}{{{\mathbb E}_p}}

\newcommand{\eps}{\epsilon}

\newcommand{\ben}{\begin{enumerate}}
\newcommand{\een}{\end{enumerate}}

\newcommand{\expect} {{\mathbb E}}

\newcommand{\Cmax} {{\Ccal}_{\rm max}}

\newcounter{countC}  
\newcounter{countR}  
\newcommand{\qn}{{\mathbb Q}_n}
\newcommand{\cn}{\Omega}

\newcommand{\R}{\Rbold}
\newcommand{\Z}{\Zbold}

\newcommand{\conn}{\leftrightarrow}

\newcommand{\bigo}{O}
\newcommand{\littleo}{o}

\title  {
        Random subgraphs of finite graphs:  \\ III.
    The phase transition for the $n$-cube
        }

\author{Christian Borgs\thanks{Microsoft Research, One Microsoft Way,
Redmond, WA 98052, USA. {\tt borgs@microsoft.com}, {\tt jchayes@microsoft.com}}
\and
Jennifer T.\ Chayes$^*$
\and
Remco van der Hofstad\thanks{Department of Mathematics and Computer Science,
Eindhoven University of Technology, P.O.\ Box  513,
5600 MB Eindhoven, The Netherlands.
{\tt rhofstad@win.tue.nl}}
\and
Gordon Slade\thanks{Department of Mathematics, University of British Columbia,
Vancouver, BC V6T 1Z2, Canada. {\tt slade@math.ubc.ca}}
\and
Joel Spencer\thanks{Department of Computer Science,
Courant Institute of Mathematical Sciences,
New York University, 251 Mercer St., New York, NY 10012, U.S.A.
{\tt spencer@cs.nyu.edu}
}}

\date\version

\begin{document}

\maketitle


\begin{abstract}
We study random subgraphs of the $n$-cube $\{0,1\}^n$, where
nearest-neighbor edges are occupied with probability $p$.
Let $p_c(n)$ be the value of $p$ for
which the expected cluster size of a fixed vertex attains the
value $\lambda 2^{n/3}$, where $\lambda$ is a small positive
constant.
Let $\epsilon=n(p-p_c(n))$.
In two previous papers, we
showed that the largest cluster inside a scaling
window given by $|\epsilon|=\Theta(2^{-n/3})$ is of size
$\Theta(2^{2n/3})$, below this scaling window it is
at most $2(\log2) n\epsilon^{-2}$, and above this scaling window
it is at most $O(\epsilon 2^n)$.
In this paper, we prove that for $p - p_c(n) \geq e^{-cn^{1/3}}$
the size of the largest cluster is at least
$\Theta(\epsilon 2^n)$, which is of the same order as the upper bound.
This provides an understanding of
the phase transition that goes far beyond that obtained by previous
authors.
The proof is based on a method that has come to be known as ``sprinkling,''
and relies heavily on the specific geometry of the $n$-cube.
\end{abstract}


\section{Introduction and results}
\label{sec-intro}

\subsection{History}
\label{sec-gt}


The study of the
random graph $G(N,p)$, defined as subgraphs of the complete graph
on $N$ vertices in which each of the possible ${N \choose 2}$ edges
is occupied with probability $p$, was initiated by
Erd\H{o}s and R\'enyi in 1960 \cite{ER60}.  They showed that
for $p = N^{-1}(1+\epsilon)$ there is a phase transition
at $\epsilon = 0$ in the sense that the size of the largest
component is $\Theta (\log N)$ for $\epsilon < 0$, $\Theta (N)$
for $\epsilon > 0$, and has the nontrivial
behavior $\Theta (N^{2/3})$ for $\epsilon = 0$.

The results of Erd\H{o}s and R\'enyi were substantially strengthened by
Bollob\'as \cite{Boll84} and {\L}uczak \cite{Lucz90}.
In particular, they showed that the model has a scaling window of
width $N^{-1/3}$, in the sense that if $p = N^{-1}(1 + \Lambda_N
N^{-1/3})$ then the size of the largest component is $\Theta
(N^{2/3})$ when $\Lambda_N$ remains uniformly bounded in $N$,
is less than $\Theta (N^{2/3})$
when $\Lambda_N \rightarrow -\infty$, and is greater than
$\Theta (N^{2/3})$ when $\Lambda_N \rightarrow +\infty$.
It is also known that inside the scaling window the
expected size of the cluster containing a given vertex is $\Theta (N^{1/3})$.

The scaling window is further characterized by the emergence of
the giant component.  When
$p=N^{-1}(1+\Lambda_N N^{-1/3})$ and $\Lambda_N\rightarrow +\infty$,
then for any constant $K>1$ the largest component will be almost
surely more than $K$ times the size of the second largest component.
When $\Lambda_N\rightarrow -\infty$ this almost surely does not happen.
However, inside the window, with $\Lambda_N=\Lambda$ fixed,
this occurs with a limiting probability strictly between
zero and one.  This is particularly striking with computer simulation.
For example, for $N=50,000$, when $\Lambda_N=-4$ the largest components are all
roughly the same size but by the ``time'' $\Lambda_N=+4$ most of them have
joined to form a dominant component several times
larger than its nearest
competitor.

In this paper,
we consider random subgraphs of the $n$-cube $\qn = \{0,1\}^n$,
where each of the nearest-neighbor edges is occupied
with probability $p$.
We
emphasize the role of the volume (number of vertices) of $\qn$ by writing
\eq
    V = |\qn| = 2^n.
\en
This model was
first
analysed in 1979
by Erd\H{o}s and Spencer \cite{ES79}, who showed that the
probability that the random subgraph is connected tends to $0$ for
$p < 1/2$, $e^{-1}$ for $p = 1/2$, and $1$ for $p> 1/2$.
More interestingly for our purposes, they showed that for
$p = n^{-1}(1 + \epsilon)$ the size of the largest component
is $o(V)$ for $\epsilon < 0$, and conjectured that it is $\Theta
(V)$ for $\epsilon > 0$.

The conjecture of Erd\H{o}s and Spencer was proved in 1982 by
Ajtai, Koml\'os and Szemer\'edi \cite{AKS82}, who thereby
established a phase transition at $\epsilon = 0$.
Their results apply for $p = n^{-1}(1+\epsilon)$ with $\epsilon$ {\em fixed}.
For $\epsilon$ fixed and negative, the largest component can be shown
by comparison with the Poisson branching process
to be a.a.s.\ of size $\bigo(n)$.
(We say that $E_n$
occurs {\it a.a.s.}\ if $\lim_{n \to \infty} \Pbold (E_n) = 1$.)
On the other hand, for $\epsilon$ fixed and positive
the largest component is at least of size $cV$ for some positive
$c = c(\epsilon)$.
To prove the latter, Ajtai, Koml\'os and Szemer\'edi introduced a method,
now known as ``sprinkling,'' which is very similar to methods
introduced at roughly the same time
in the context of percolation on ${{\mathbb Z}^d}$
by Aizenman, Chayes, Chayes, Fr\"ohlich and Russo \cite{ACCFR83}. We
will use a variant of sprinkling in this paper.
Very recently, Alon, Benjamini and Stacey \cite{ABS02} used the
sprinkling technique to extend the Ajtai, Koml\'os and
Szemer\'edi result \cite{AKS82} to
subgraphs of transitive finite graphs of
high girth.  These results, while applicable to much more than the
$n$-cube, also hold only for $\epsilon$ fixed.

A decade
after the Ajtai, Koml\'os and Szemer\'edi work, Bollob\'as,
Kohayakawa and \L uczak \cite{BKL92}
substantially refined
their
result,
in particular studying the behavior of the largest cluster as
$\epsilon \rightarrow 0$.
Let $|C(x)|$ denote the size of the cluster of $x$,
let $\Cmax$ denote a cluster of maximal size, and let
    \eq\lbeq{Cmaxdef}
    |\Cmax| = \max\{|C(x)| :  x \in \qn \}
    \en
denote the maximal cluster size.
We again take
$p= n^{-1}(1+\epsilon)$, and now assume that $\epsilon \to 0$ as $n \to\infty$.
In \cite[Corollary~16, Theorem~28]{BKL92}
(with somewhat different notation), it is proved that
for $\epsilon \leq -(\log n)^2(\log \log n)^{-1}n^{-1/2}$,
    \eq\lbeq{BKL-e<0}
    |\Cmax| = \frac{2 \log V}{\epsilon^2} \big(1 + o(1)\big)
    \quad a.a.s. \;\;\mbox{as $n \to \infty$},
    \en
and that for $\epsilon \geq 60 (\log n)^3 n^{-1}$,
    \eq\lbeq{BKL-e>0}
    |\Cmax| = 2\epsilon V \big(1 + o(1)\big)
    \quad a.a.s. \;\;\mbox{as $n \to \infty$}.
    \en
Thus, $\epsilon \geq 60 (\log n)^3 n^{-1}$ is supercritical.
In addition, it is shown in \cite[Theorem~9]{BKL92} that
the right side of \refeq{BKL-e<0} is an upper bound on $|\Cmax|$
provided that $p< (n-1)^{-1} - e^{-o(n)}$, and hence such $p$
are subcritical.

In recent work \cite{BCHSS04a,BCHSS04b}, we developed a general theory
of percolation on
connected transitive finite graphs that applies to $\qn$ and
to various high-dimensional tori.  The theory is based on the view that
the phase transition on many high-dimensional graphs should have similar
features to the phase transition on the complete graph.  In particular,
the largest cluster should have size $\Theta(V^{2/3})$ in a scaling window
of width $\Theta(V^{-1/3})$, and should have size $o(V^{2/3})$ below the
window and  size $\Theta(V)$ above the window.
Note that the bounds of \cite{BKL92},
while much sharper than those established in \cite{AKS82},  are
still far from establishing this behavior.

We will review the
results of \cite{BCHSS04a,BCHSS04b} in detail below,
as they apply to $\qn$.
These results do not give a lower bound on the largest cluster above
the scaling window, and our primary purpose in this
paper is to provide such a bound.  We define a critical threshold
$p_c(n)$ and prove that the largest cluster has size
$\Theta([p-p_c(n)]nV)$ for $p-p_c(n) \geq e^{-cn^{1/3}}$.
This falls short of proving a bound for all $p$ above a window of
width $V^{-1/3}$, but it greatly extends the range of $p$
covered by the Bollob\'as, Kohayakawa and \L uczak bound \refeq{BKL-e>0}.

\subsection{The critical threshold}

The starting point in \cite{BCHSS04a}
is to define the critical threshold in terms of
the {\em susceptibility}\/ $\chi(p)$, which is defined to be
the expected size of the cluster of a given vertex:
    \eq\lbeq{chidef}
    \chi(p) ={\mathbb E}_p|C(0)|.
    \en
For percolation on $\Zd$, $\chi(p)$ diverges to infinity as $p$ approaches
the critical point from below.  On $\qn$,
the function $\chi$ is strictly monotone increasing on the interval
$[0,1]$, with $\chi(0)=1$ and $\chi(1)=V$.  In particular, $\chi(p)$ is
finite for all $p$.

For $G(N,p)$, the susceptibility is $\Theta(N^{1/3})$ in the scaling
window.  For $\qn$, the role of $N$ is played by
$V=2^n$, so we could expect by analogy that $p_c(n)$ for the $n$-cube
should be roughly equal to the $p$ that solves
$\chi(p) = V^{1/3} = 2^{n/3}$.
In \cite{BCHSS04a}, we defined the critical threshold
$p_c = p_c(n) = p_c(n;\lambda)$ by
\eq
\lbeq{pcdef}
     \chi(p_c) = \lambda V^{1/3},
\en
where $\lambda$ is a small positive constant.
The flexibility in the choice of $\lambda$ in
\refeq{pcdef} is connected with the fact that the phase transition in a finite
system is smeared over an interval rather than occurring at a sharply
defined threshold, and any value in the transition interval could be
chosen as a threshold.  Our results
show that $p_c$ defined by \refeq{pcdef}
really is a critical threshold for percolation on $\qn$.

The triangle condition plays an important role in the analysis of
percolation on $\Zd$ for large $d$ \cite{AN84,BA91,HS90a}, as well
as on infinite non-amenable graphs \cite{Scho01}.
For $x,y \in \qn$, let $\{x\conn y\}$ denote the event that $x$ and
$y$ are in the same cluster, and let $\tau_p(x,y) = \Pbold_p(x \conn y)$.
The {\em triangle diagram} is defined by
\eq
    \nabla_p(x,y) = \sum_{w,z \in \qn} \tau_p(x,w)\tau_p(w,z)\tau_p(z,y).
\en
The {\em triangle condition} is the statement that
\eq
\lbeq{tc}
    \max_{x,y \in \qn}\nabla_{p_c(n)}(x,y) \leq \delta_{x,y} +a_0
\en
where $a_0$ is less than
a sufficiently small constant.  The {\em stronger triangle condition}
is the statement that there are positive constants
$K_1$ and $K_2$ such that $\nabla_p(x,y) \leq \delta_{x,y} + a_0$
uniformly in $p \leq p_c(n)$, with
\eq
\lbeq{a0}
    a_0 = a_0(p) = K_1 n^{-1} + K_2 \chi^3(p) V^{-1}.
\en
If we choose $\lambda$ sufficiently small
and $n$ sufficiently large,
then the stronger triangle condition implies the triangle condition.

In \cite{BCHSS04b}, we proved the stronger triangle condition for $\qn$.
In \cite[Theorems~1.1--1.5]{BCHSS04a},
we derived several consequences of the stronger triangle
condition in the context of general finite connected transitive graphs.
For the $n$-cube, these results (see \cite[Theorems~1.1, 1.5]{BCHSS04a}
imply that there is a $\lambda_0 >0$
and a $b_0>0$ (depending on $\lambda_0$)
such that if $0<\lambda \leq \lambda_0$ then
    \eq\lbeq{omegapc-1}
    1- \lambda^{-1} V^{-1/3}
    \leq n p_c(n) \leq
    1 + b_0n^{-1}.
    \en
%
In addition, given $0<\lambda_1 <\lambda$, let $p_1$ be defined by
$\chi(p_1)=\lambda_1 V^{1/3}$.  Then
\eq
\lbeq{p2p1}
    \frac{\lambda -\lambda_1}{\lambda_1\lambda } \frac{1}{V^{1/3}}
    \leq
    n(p_c(n)-p_1 )
    \leq
    \frac{\lambda -\lambda_1}
    {\lambda_1\lambda } \frac{1}{V^{1/3}}[1+b_0n^{-1}].
\en
Thus decreasing $\lambda$ to $\lambda_1$
shifts $p_c(n)$ only by $O(n^{-1} V^{-1/3})$, so $p_1$ remains in the
scaling window.

The asymptotic formula $p_c(n) = n^{-1} + \bigo(n^{-2})$ of
\refeq{omegapc-1} is improved in \cite{HS03a} to
        \eq
        \lbeq{pasy}
        p_c(n) = \frac{1}{n}+\frac 1{n^2} + \frac{7}{2n^3}
        +\bigo(n^{-4}).
        \en
Presumably there is an asymptotic expansion to all orders, so
that there are real numbers
$a_i$ ($i \geq 1$) such that for each $s\geq 1$
        \eq
        \lbeq{pcasyexp}
        p_c(n) = \sum_{i=1}^s a_i n^{-i} + \bigo(n^{-(s+1)}).
        \en
Except in the unlikely event that $a_i=0$ for all $i$ sufficiently
large, so that the expansion is actually a finite polynomial in $n^{-1}$,
we see from \refeq{pcasyexp} that for every $s \geq 1$,
the truncated expansion $\sum_{i=1}^s a_i n^{-i}$ lies outside an
interval of width
$V^{-1/3}=2^{-n/3}$ centered at $p_c(n)$, for large $n$.
The non-perturbative definition
\refeq{pcdef} of $p_c(n)$ therefore tracks the scaling window
more accurately than any polynomial in $n^{-1}$ can ever do.  Moreover,
as is discussed in more detail in \cite{HS03a},
we expect that the full expansion $\sum_{i=1}^\infty a_i n^{-i}$ is
a divergent series.  This would mean that, given $n$, if we
take $s$ large depending on $n$ then the truncated series
$\sum_{i=1}^s a_i n^{-i}$ would be meaningless---possibly not even
lying in the interval $[0,1]$.

Bollob\'as, Kohayakawa and \L uczak raised the question
of whether the
critical value might be {\em equal}\/ to $(n-1)^{-1}$.  Note that $n-1$ is
the forward branching ratio of the ``tree approximation'' to the
$n$-cube, so that this suggestion would mean that the tree
approximation would give the correct critical value. Note also
that, because of the smearing of the random graph critical point
by the scaling window of width $N^{-1/3}$, on the random graph
there is no distinction between the critical values $p_c(N) =
(N-1)^{-1}$ and $p_c(N) = N^{-1}$, so that the the tree
approximation does give a correct critical value in that case.
However, for the $n$-cube, our picture that
the width of the scaling window is $\Theta(V^{-1/3})=\Theta(2^{-n/3})$
implies that there
is a real distinction between the values $n^{-1}$ and
$(n-1)^{-1}$, and \refeq{pasy}
implies that both lie {\em below}\/ the critical window.

\subsection{In and around the scaling window}

Given $p \in [0,1]$, let
$\epsilon =\epsilon(p) \in \R$ be defined by
\eq
\lbeq{epdef}
    p = p_c(n) + \frac{\epsilon}{n}.
\en
We say that $p$ is {\em below} the window (subcritical) if
$\epsilon V^{1/3} \to -\infty$, {\em above} the window
(supercritical) if $\epsilon V^{1/3} \to
\infty$, and {\em inside} the window if
$|\epsilon| V^{1/3}$ is uniformly bounded
in $n$.  In this section, we summarize and rephrase the results stated in
\cite[Theorems~1.2--1.5]{BCHSS04a}, as they apply to $\qn$.
We have stated slightly weaker results than those obtained in
\cite{BCHSS04a}, in order to simplify the statements.
These results where proved in
\cite{BCHSS04a} assuming the triangle condition (or the stronger triangle
condition for \cite[Theorem~1.5]{BCHSS04a}, and the stronger triangle
condition was established in \cite{BCHSS04b}.

\begin{theorem}[Below the window]
\label{main-thm-sub} Let $\lambda \leq \lambda_0$ and
$p=p_c(n) -\epsilon n^{-1}$ with $\epsilon \geq 0$.
Then
    \eq\lbeq{chiasy}
    \frac{1}{\lambda^{-1} V^{-1/3}+\epsilon}
    \leq
    \chi(p)
    \leq \frac{1}%
    {\lambda^{-1} V^{-1/3}+[1- a_0] \epsilon  },
    \en
with $a_0=a_0(p)$ given by \refeq{a0}.
Moreover, if $\epsilon V^{1/3} \to \infty$ as $n \to \infty$, then
\eq
    \chi(p) = \frac{1}{\epsilon} [1+o(1)],
\en
    \eq\lbeq{cmaxbd1}
    \frac{1}{10^4 \epsilon} [1+o(1)]
    \leq
    \Exp\big(|\Cmax|\big)
    \leq
    \frac{2\log V}{\epsilon^2}[1+o(1)],
    \en
\eq\lbeq{cmaxbd2}
    \frac{1}{3600 \epsilon^2} [1+o(1)] \leq
    |\Cmax|\leq \frac{2\log V}{\epsilon^2}[1+o(1)]
    \quad \text{a.a.s.}
\en
\end{theorem}

In the language of critical exponents, the above bounds on $\chi(p)$
correspond to $\gamma =1$.
Since $p_c(n) >
(n-1)^{-1}$ for large $n$ by \refeq{pasy}, the upper bound of \refeq{cmaxbd1}
extends the range of $p$ covered by the Bollob\'as, Kohayakawa and \L uczak
upper bound of \refeq{BKL-e<0}
from $p < (n-1)^{-1}-e^{-\littleo(n)}$ to all $p$ below the window.
We conjecture that the upper bound of \refeq{cmaxbd1}
is actually sharp for all
$p$ that are not exponentially close to $p_c(n)$.
This was proved in \cite{BKL92} (see \refeq{BKL-e<0}) for
$\epsilon \geq (\log n)^2(\log \log n)^{-1}n^{-1}$.
We also conjecture that this behavior can be extended appropriately
to cover a larger range of $\epsilon$, as follows.

    \begin{conj}
    \label{conj-LCsub} Let $p = p_c(n) - \epsilon n^{-1}$
    with $\lim_{n \to \infty} \epsilon = 0$ and
    $\lim_{n \to \infty} \epsilon e^{\delta n} = \infty$ for
    every $\delta > 0$.  Then
        \eq
        \lbeq{LCasy}
        |\Cmax|= \frac{2 \log{V} }{\epsilon^2}[1+\littleo(1)]
        \qquad  \text{ a.a.s. }
        \en
    If we assume instead that $\lim_{n \to \infty} \epsilon = 0$ and
    $\lim_{n \to \infty} \epsilon 2^{n/3} = \infty$, then
        \eq
        \lbeq{LCasy2}
        |\Cmax|= \Theta\big(
    \frac{2\log(\epsilon^3 V)}{\epsilon^2}\big)
        \qquad \text{ a.a.s. }
        \en
    \end{conj}

Note that \refeq{LCasy2}
reduces to \refeq{LCasy} when $\epsilon$ is not exponentially small.
The asymptotic behavior \refeq{LCasy2}, with $V$ replaced by $N$,
is known to apply to
the random graph for $N^{-1/3} \ll N(p_c-p) \ll 1$;
see \cite[Theorem~5.6]{JLR00}.

For $k \geq 0$, let
\eq
    P_{\geq k} = \Pbold_p(|C(0)| \geq k).
\en

\begin{theorem}[Inside the window]
\label{main-thm-critical}

Let $\lambda \leq \lambda_0$ and $\Lambda<\infty$.
Let $p=p_c +\cn^{-1}\epsilon$
with $|\epsilon|\leq\Lambda V^{-1/3}$.
There are finite
positive constants $b_1,\dots,b_8$ such that the
following statements hold.

\noindent
i) If $k\leq b_1 V^{2/3}$, then
\eq\lbeq{clszdis}
\frac{b_2}{\sqrt k}
\leq
P_{\geq k}(p)
\leq
\frac{b_3}{\sqrt k}.
\en

\noindent
ii) 
\eq\lbeq{LCEBd1-win}
    {b_4}V^{2/3}
    \leq
    \Exp\big(|\Cmax|\big)
    \leq
    {b_5}V^{2/3}
    \en
and, if $\omega\geq 1$, then
    \eq\lbeq{LCBd1-win}
    \Pro\Big(
    \omega^{-1} V^{2/3}\leq |\Cmax|\leq \omega V^{2/3}
        \Big)
    \geq 1-\frac{b_6}\omega.
    \en

\noindent
iii)
    \eq\lbeq{chiasy-win}
    b_7 V^{1/3}
    \leq \chi(p)\leq
    b_8 V^{1/3}.
    \en

\noindent In the above statements, the constants $b_2$ and $b_3$
can be chosen
to be independent of $\lambda$ and $\Lambda$, the constants
$b_5$ and $b_8$ depend on $\Lambda$ and not $\lambda$, and the constants
$b_1$, $b_4$, $b_6$ and $b_7$ depend on both
$\lambda$ and $\Lambda$.

\end{theorem}

In terms of critical exponents,
\refeq{clszdis} says that $\delta = 2$.

    \begin{theorem} [Above the window]
    \label{thm-supercritical}
Let $\lambda \leq \lambda_0$ and $p = p_c + \epsilon n^{-1}$ with
$\epsilon V^{1/3} \to \infty$.
Then
    \eq
    \lbeq{chiasysup}
    \chi(p)
    \leq 162 \epsilon^2 V,
    \en
\eq
    \lbeq{bound1-cmax}
        \Exp(|\Cmax|)
        \leq  28\epsilon V,
        \en
and, for all $\omega >0$,
\eq\lbeq{cmax.2A}
    \Pro\Big(|\Cmax|\geq \omega \epsilon V\Big)
    \leq \frac{\rm{const}}{\omega}.
    \en
\end{theorem}

A refinement of \refeq{cmax.2A} will be given in Section~\ref{sec-giantub}.

To see that there is a phase transition at $p_c(n)$, we need an upper bound
on the maximal cluster size in the subcritical phase and a lower bound
in the supercritical phase.  The former is given in Theorem~\ref{main-thm-sub}
but the latter is not part of Theorem~\ref{thm-supercritical}.

\subsection{Main result}
\label{sec-mr}

Our main result is the following theorem, which
is proved in
Section~\ref{sec-sprinkling}.
Theorem~\ref{thm-mr} provides the missing lower
bound for $\epsilon \geq e^{-cn^{1/3}}$.  This restriction on
$\epsilon$ is an artifact of our proof and we believe the theorem remains
valid as long as $\epsilon V^{1/3}\to\infty$; see Conjecture~\ref{conj-2}.
To fully establish the picture that there is a scaling window of width
$\Theta(V^{-1/3})$, it would be necessary to extend Theorem~\ref{thm-mr}
to cover this larger range of $\epsilon$.

    \begin{theorem}
    \label{thm-mr}
There are $c,c_1>0$ and $\lambda_0 >0$ such that the following hold for
all $0<\lambda \leq \lambda_0$ and all
$p = p_c + \epsilon n^{-1}$ with $e^{-cn^{1/3}} \leq \epsilon \leq 1$:
    \eqalign
    \lbeq{mr1}
    |\Cmax| & \geq c_1 \epsilon 2^n
    \quad\text{ a.a.s. as }n\rightarrow \infty,
    \\
    \lbeq{mr2}
    \chi(p)
    & \geq [1+\littleo(1)](c_1 \epsilon)^2 2^n
    \quad\text{ as }n\rightarrow \infty.
    \enalign
\end{theorem}

It is interesting to examine the
approach to the critical point with $|p - p_c(n) |$ of order $n^{-s}$
for different values of $s$.
Our results give a hierarchy of bounds as $s$ is varied.
For example, it follows from Theorem~\ref{main-thm-sub} that
for $p=p_c(n) - \delta n^{-s}$ with $s>0$ and $\delta >0$,
   \eqalign
    \lbeq{chiasy2}
        \chi(p) & = n^{s-1} \delta^{-1}[1+\bigo(n^{-1})].
    \enalign
On the other hand, for $p=p_c(n) + \delta n^{-s}$ with $s>0$ and $\delta >0$,
it follows from Theorems~\ref{thm-supercritical} and \ref{thm-mr} that
    \eqalign
    \lbeq{chiasysup2}
    \chi(p)
    = \Theta( \delta^{2} n^{2(1-s)} 2^n).
    \enalign
Related bounds follow for $|\Cmax|$.
Thus there is a phase transition on scale
$n^{-s}$ for any $s \geq 1$.

\subsection{More conjectures}

For $\epsilon \geq 60(\log n)^3 n^{-1}$, Bollob\'as, Kohayakawa
and \L uczak \cite{BKL92}
proved (see \refeq{BKL-e>0}) that $|\Cmax| = 2\epsilon 2^n(1+o(1))$ a.a.s.
We conjecture that this
formula holds for all $\epsilon$ above the window.
This behavior has been proven for the random graph above the scaling window;
see \cite[Theorem~5.12]{JLR00}.

\begin{conj}
\label{conj-2}
Let $p=p_c(n)+ \epsilon n^{-1}$ with
$\epsilon >0$, $\lim_{n \to \infty} \epsilon = 0$
and $\lim_{n \to \infty} \epsilon V^{1/3}=\infty$.  Then
    \eqalign
    \lbeq{LCbdsup3}
    |\Cmax| &= 2 \epsilon  V[1+o(1)]\qquad
    \text{ a.a.s.},
    \\
    \lbeq{chiasysup3}
    \chi(p) &= 4\epsilon^2 V[1+o(1)].
    \enalign
\end{conj}

The constants in the above conjecture
can be motivated by analogy to the Poisson branching process
with mean $\lambda$.  There the critical point is $\lambda=1$ whereas
we have a critical point $p_c(n)$.  An increase in $p$ beyond $p_c(n)$ by
$\epsilon n^{-1}$ increases the average number of neighbors
of a vertex by $\epsilon$, which we believe corresponds to the Poisson
branching process with mean $1+\epsilon$.  This process is infinite
with probability $\sim 2\epsilon$.  When we generate the
component of $x$ in $\qn$ it cannot, of course, be infinite,
but with probability $\sim 2\epsilon$ it will not die quickly.
Consider components of $x,y\in \qn$ that do not die quickly.
We believe that these components will not avoid each other.
Rather, all of them will coalesce to form a component $\Cmax$
of size $2\eps V$.  Finally,
with probability $\sim 2\epsilon$ a given vertex $0$ lies in $\Cmax$
and this contributes $\sim (2\epsilon)(2\epsilon V)$ to $\chi(p)$,
which we believe is the dominant contribution.
Note also that for any positive {\em constant}\/ $\epsilon$,
it is shown in \cite[Theorem~29]{BKL92}
that $|\Cmax|\sim
a N$ where $a$ is the probability that the Poisson branching
process with mean $1+\epsilon$ is infinite.

For the random graph $G(N,p)$,
outside the scaling window there is an intriguing {\em duality}
between the subcritical and supercritical phases
(see \cite[Section~10.5]{AS00}, and, for a more general setting, see
\cite{MR98}).
Let $p= N^{-1}(1+\epsilon)$
lie above the scaling window for $G(N,p)$, and to
avoid unimportant issues, assume that $\epsilon =o(1)$.
Almost surely, there is a dominant component
of size $\sim 2\epsilon N$.
Remove this component from the graph, giving $G^-$.
Then $G^-$ behaves like the random graph in the subcritical
phase with probability $p'= p_c(1-\epsilon)$.
In particular, the size of the largest component of $G^-$ (the second largest
component of $G$) is given asymptotically by the size of the
largest component for $p'$.

We believe that this duality holds for random subgraphs of
$\qn$ as well.  Let $\Ccal_2$ denote the second largest cluster, and
set $p=n^{-1}(1+\epsilon)$.  Bollob\'as, Kohayakawa and \L uczak
\cite{BKL92} showed that if $\epsilon\rightarrow 0$ and
$\epsilon \geq 60 (\log n)^3 n^{-1}$ then
$|\Ccal_2|\sim (2\log V)\epsilon^{-2}$
a.a.s.  Note that this matches the behavior of
$|\Cmax|$ for $p= p_c(n) -\epsilon n^{-1}$ given in
Conjecture~\ref{conj-LCsub}.  The following conjecture, which we
are far from able to show using our present methods even for
$\epsilon \geq e^{-cn^{1/3}}$, claims an extension of the result
of \cite{BKL92} to all $p$ above the window.

\begin{conj}
\label{conj-C2}
Let $p=p_c(n)+ \epsilon n^{-1}$ with
$\epsilon >0$, $\lim_{n \to \infty} \epsilon = 0$
and $\lim_{n \to \infty} \epsilon V^{1/3}=\infty$.  Then as $n \to \infty$,
the size of the second largest cluster $\Ccal_2$
is
\eq
    |\Ccal_2|  = \frac{2\log V}{\epsilon^2}[1+o(1)]
    \quad \mbox{ a.a.s. }
\en
\end{conj}

Finally, we consider the largest cluster inside the scaling window.
For $G(N,p)$, it is known that the largest cluster inside the
window has size $X N^{2/3}$ where $X$ is a positive random variable
with a particular distribution.  Similarly, we expect that
for $\qn$ the largest
cluster inside the window has size $YV^{2/3}$ for some positive
random variable $Y$.  Our current methods are not sufficient
to prove this.

\section{Proof of Theorem~\ref{thm-mr}}
\label{sec-sprinkling}

In this section, we prove Theorem~\ref{thm-mr}
by showing that there is a $c_1>0$ such that
when $e^{-cn^{1/3}} \leq \epsilon \leq 1$,
    \eq
    \lbeq{LCbdsupz}
    |\Cmax| \geq c_1 \epsilon 2^n
    \quad\text{ a.a.s. }
    \en
and
    \eq\lbeq{chiasysupz}
    \chi(p) \geq (c_1 \epsilon)^2 2^n [1+\littleo(1)].
    \en
The proof of \refeq{LCbdsupz} is based on
the method of sprinkling and is given in
Sections~\ref{sec-giantub}--\ref{sec-giantlb}.
The bound \refeq{chiasysupz} is an elementary
consequence and is given in Section~\ref{sec-chibdsuper}.

\subsection{The percolation probability and sprinkling}
\label{sec-giantub}

For percolation on $\Z^d$, the value of $p$
for which $\chi(p)$ becomes infinite is the same as
the value of $p$ where the percolation probability ${\mathbb
P_p}(|C(0)|=\infty)$ becomes positive \cite{AB87,Mens86}.
For $\qn$, there can be no
infinite cluster, and the definition of the percolation probability
must be modified.  For $p=p_c(n)+\epsilon n^{-1}$ with
$\epsilon >0$,
we defined the percolation probability in \cite{BCHSS04a} by
    \eq\lbeq{theta_a}
    \theta_\alpha(p) =
    \Pro(|C(0)| \geq N_\alpha)
    =
    P_{\geq N_\alpha},
    \en
where
    \eq\lbeq{N_a}
    N_\alpha=  N_\alpha(p) =
    \frac 1{\epsilon^2}
    \big(\epsilon V^{1/3}\big)^\alpha
    \en
and $\alpha$ is a fixed parameter in $(0,1)$.
The definition \refeq{theta_a} is motivated as follows.
According to Conjectures~\ref{conj-2}--\ref{conj-C2},
above the window
the largest cluster has size $|\Cmax| =
2 \epsilon  V[1+o(1)]$ a.a.s., while the second largest has size
$|\Ccal_2|=2\epsilon^{-2}\log V[1+o(1)]$.
According to this, above the window
$|\Ccal_2| \ll N_\alpha \ll
|\Cmax|$, so that a cluster of size at least $N_\alpha$ should
in fact be maximal, and $\theta_\alpha(p)$ should correspond to the
probability that the origin is in the maximal cluster.
(The above reasoning suggests
the range $0<\alpha<3$ rather than $0<\alpha <1$, but the analysis
of \cite{BCHSS04a} requires the latter restriction.)

Let $0<\alpha<1$.  The combination of \cite[Theorem~1.6]{BCHSS04a} with
the verification of the triangle condition in \cite{BCHSS04b}
implies that there are
positive constants $b_9$, $b_{10}$
such that
        \eq\lbeq{bound1}
         b_{10}\epsilon
         \leq
        \theta_\alpha (p)
        \leq
        27\epsilon,
        \en
where the lower bound holds when $b_9 V^{-1/3}\leq \epsilon\leq 1$
and the upper bound holds when $\epsilon \geq V^{-1/3}$.
In addition, there are  positive $b_{11}$, $b_{12}$
such that if $\max\{b_{12}V^{-1/3}, V^{-\eta}\}\leq
\epsilon\leq 1$, where
$\eta=\frac 13\frac{3-2\alpha}{5-2\alpha}$,
then
    \eq
    \lbeq{LCbdsup}
    \Pro\Big(
    |\Cmax|\leq [1 + (\epsilon V^\eta)^{-1}]\theta_\alpha(p) V
    \Big)
    \geq 1- \frac{b_{11}}{(\epsilon V^\eta)^{3-2\alpha}}.
    \en
In the above statements, the constants  $b_9$,
$b_{10}$, $b_{11}$ and $b_{12}$ depend on both $\alpha$ and
$\lambda$.
Note that although \refeq{LCbdsup} does not obtain the precise constant
of \refeq{BKL-e>0} found by Bollob\'as, Kohayakawa
and \L uczak, it does extend the range of $p$ from
$p \geq (n-1)^{-1} + 60 (\log^2 n)n^{-2}$
to $p \geq p_c(n) + 2^{-\eta'}$ for any $\eta' < \eta$.
Also, note that the combination of \refeq{LCbdsup} and \refeq{bound1}
gives a refinement of \refeq{cmax.2A}.

Let
\lbeq{Zdef}
    \eq
    Z_{\geq N_\alpha} = \sum_{x\in \qn} I[|C(x)|\geq N_\alpha]
    \en
denote the number of vertices in ``moderately'' large components.  Then
$\expect_p(Z_{\geq N_\alpha})=\theta_\alpha(p) V$ and hence,
by \refeq{bound1}, above the window
$\expect_p(Z_{\geq N_\alpha})= \Theta(\epsilon V)$.
In the proof of \cite[Theorem~1.6~ii)]{BCHSS04a}, it is shown that
\eq
\lbeq{v1}
    \Pbold_p \big(|Z_{\geq N_\alpha}- V\theta_\alpha(p)|
    \geq (\epsilon V^\eta)^{-1} V\theta_\alpha(p) \big)
    \leq
    \frac{b_{11}}{(\epsilon V^\eta)^{3-2\alpha}}
\en
for percolation on an arbitrary finite
connected transitive graph that obeys the
triangle condition, and hence \refeq{v1} holds for $\qn$.
For $\epsilon \geq e^{-cn^{1/3}}$ and fixed $\alpha \in (0,1)$, there are
therefore positive constants $\eta_1, \eta_2, A$
such that
\eq
\lbeq{varest2}
    \Pbold_p\big( |Z_{\geq N_\alpha}-V\theta_\alpha(p)|
    \geq V^{1-\eta_1} \theta_\alpha(p) \big)
    \leq
    A V^{-\eta_2}.
\en
This shows that $Z_{\geq N_\alpha}$ is typically close to
its expected value $V\theta_\alpha(p)$.
For $V$ sufficiently large and $e^{-cn^{1/3}}
\leq \epsilon \leq 1$, it follows from \refeq{varest2} and the lower bound
of \refeq{bound1} that
\eq
\lbeq{varest}
    \Pbold_p\big( Z_{\geq N_\alpha} \geq
    2 b_{10}\epsilon V \big)
    \leq
    A V^{-\eta_2}.
\en

The proof of \refeq{LCbdsupz} is based on the following sketch.
Let $p \geq p_c(n)+\epsilon n^{-1}$ with $e^{-cn^{1/3}}\leq \epsilon \leq 1$.
Let $p^- = \frac 12 (p-p_c(n))$, and define $p^+$ by
$p^- + p^+ - p^-p^+ = p$.  Then a percolation configuration with bond
density $p$ can be regarded as the union of two
independent percolation configurations having bond densities $p^-$ and
$p^+$.  The additional bonds due to the latter are regarded
as having been ``sprinkled'' onto the former.
For percolation with bond density $p^-$, it follows from \refeq{varest}
that a positive fraction of the vertices lie in moderately large components.
We then use the specific geometry of $\qn$, in a crucial way,
to argue that
after a small sprinkling of additional bonds a positive fraction
of these vertices will be joined together
into a single giant component,
no matter how the vertices in the large components are arranged.
Our restriction $\epsilon \geq e^{-cn^{1/3}}$ enters in this last step.

\subsection{The largest cluster}
\label{sec-giantlb}

In Proposition~\ref{sprinkle} below, we
prove  the lower bound on $|\Cmax|$
of \refeq{LCbdsupz}.  In preparation
for Proposition~\ref{sprinkle}, we state four lemmas.  The first lemma
uses a very special geometric property of
$\qn$.  For its statement, given
any $X\subseteq \qn$ and positive integer $d$, we denote
the ball around $X$ of radius $d$ by
\eq
    B[X,d]= \{y\in \qn: \exists x\in X \mbox{ such that } \rho(x,y) \leq d\},
\en
where $\rho(x,y)$ denotes the graph distance between $x$ and $y$.

\begin{lemma}[Isoperimetric Inequality] \label{isoperimetric}
If $X\subseteq \qn$ and $|X|\geq \sum_{i\leq u}{n\choose i}$
then
\eq  |B[X,d]| \geq \sum_{i\leq u+d}{n\choose i} . \en
\end{lemma}

Lemma~\ref{isoperimetric} is proved in Harper \cite{Harp66}.  Bollob\'as
\cite{Boll86} is a very readable and more modern reference.
The result of Lemma~\ref{isoperimetric}
may be seen to be best possible by taking $X=B[\{v\},d']$
for any fixed $v\in \qn$, so that $B[X,d]=B[\{v\},d'+d]$.  For
asymptotic calculations we use the inequality of the following lemma.

\begin{lemma}[Large Deviation] \label{largedeviation}
For $\Delta >0$,
\eq \sum_{i\leq \frac{n-\Delta}{2}} {n\choose i} =
\sum_{i\geq \frac{n+\Delta}{2}} {n\choose i} \leq 2^ne^{-\Delta^2/2n}.
\en
\end{lemma}

\proof
The first two terms are equal by the symmetry of Pascal's triangle.
Dividing by $2^n$, the inequality may be regarded as the large
deviation inequality
\eq
    \Pr[S_n\geq \Delta] \leq e^{-\Delta^2/2n},
\en
where $S_n=\sum_{i=1}^n X_i$ with the $X_i$ independent
random variables with $\Pbold (X_i = \pm 1)=\frac 12$.
A simple proof
of this basic inequality is given in \cite[Theorem~A.1.1]{AS00}.
\qed

\begin{lemma}[Big Overlap]\label{bigintersection}
Let $\Delta,\epsilon >0$ satisfy $e^{-\Delta^2/2n}<\frac{\epsilon}{2}$.
Let $S,T\subseteq \qn$ with $|S|,|T|\geq \epsilon 2^n$.  Then
\eq |B[S,\Delta]\cap T| \geq \frac{1}{2} |T| . \en
\end{lemma}

\proof
From Lemma~\ref{largedeviation}, $|S|\geq \sum_{i\leq (n-\Delta)/2}
{n\choose i}$.  Hence, by Lemma~\ref{isoperimetric},
$|B[S,\Delta]| \geq \sum_{i\leq (n+\Delta)/2}{n\choose i}$ (intuitively, we
have crossed the equator).  Therefore, by Lemma~\ref{largedeviation},
\eq |\qn \setminus B[S,\Delta]| \leq \sum_{i< \frac{n-\Delta}{2}} {n \choose i}
<\frac{\epsilon}{2}
2^n\leq \frac{1}{2}|T|,  \en
and so $B[S,\Delta]$ must overlap at least half of $T$.
\qed

\begin{lemma}[Many Paths]\label{manypaths}
Let $\Delta,\epsilon >0$ satisfy $e^{-\Delta^2/2n}<\frac{\epsilon}{2}$.
Let $S,T\subseteq \qn$ with $|S|,|T|\geq \epsilon 2^n$.  Then
there is a collection of $\frac{1}{2}\epsilon 2^nn^{-2\Delta}$ vertex
disjoint paths from $S$ to $T$, each of length at most $\Delta$.
\end{lemma}

\proof
Set $T_1=B[S,\Delta]\cap T$.
By Lemma~\ref{bigintersection}, $|T_1|\geq \frac{1}{2}\epsilon 2^n$.
Let $T_2\subseteq T_1$ be a maximal subset such that no $x,y\in T_2$
are within distance $2\Delta$ of each other.
Every $y\in T_1$ must lie in a ball
of radius $2\Delta$ around some $x\in T_2$ and each such ball has
size at most $n^{2\Delta}$ (using a crude upper bound), so
$|T_2|\geq n^{-2\Delta}|T_1|$.  For each $x\in T_2$ there is a path
of length at most $\Delta$ to some $z\in S$ and the paths from
$x,y\in T_2$ must be disjoint as otherwise $x,y$ would be at distance
at most $2\Delta$.
\qed

\medskip
Now we use Lemma~\ref{manypaths}
and \refeq{varest} to prove \refeq{LCbdsupz}.

\begin{prop}[Sprinkling]
\label{sprinkle}
There are absolute positive constants $c_1, \beta$ such that
\eq
\lbeq{Cmaxlb}
    \Pbold ( |\Cmax| \leq c_1 \epsilon 2^n)
    \leq
    2^{-\beta n}.
\en
whenever $e^{-cn^{1/3}} \leq \epsilon \leq 1$.
In particular, $  |\Cmax| \geq c_1\epsilon 2^n$ a.a.s.

\end{prop}

\begin{proof}  As usual, we write $p=p_c(n)+\epsilon n^{-1}$.
Let $p^-$ be such that
\eq
    p^- + \frac{\epsilon}{2n} - \frac{\epsilon}{2n}p^- = p .
\en
Note that $p^-=p_c(n)+ \epsilon/(2n) + o(\epsilon/n)$. We consider
the random subgraph of $\qn$ with probability $p$ as the union
of the random subgraph $G^-$ with probability $p^-$ and the
random subgraph $H$ (the sprinkling) with probability $\epsilon/(2 n)$.
Crucially,
$G^-$ and $H$ are chosen independently.

Let $C_i$ ($i\in I$)
denote the components of $G^-$ of size at least $2^{\alpha n/3}$.
Set
\eq
    D = \bigcup_{i\in I}C_i \quad \mbox{and}
    \quad M = |D| ,
\en
and note that
\eq
\lbeq{Nbd}
    N_\alpha = \frac{1}{\epsilon^{2-\alpha}} V^{\alpha/3}
    \geq V^{\alpha /3} = 2^{\alpha n/3}.
\en
Since $M \geq Z_{\geq N_\alpha}$ by \refeq{Nbd}, it follows from
\refeq{varest} that there is an absolute positive
constant $c_2$ such that
\eq
    \Pbold_{p^-}(M\geq c_2\epsilon 2^n)
    \geq
    \Pbold_{p^-}(Z_{\geq N_\alpha} \geq c_2\epsilon 2^n) \to 1
\en
exponentially rapidly in $n$.  Thus we may assume that
$G^-$ has $M\geq c_2\epsilon 2^n$.
It suffices to show that the probability that at least $\frac{1}{3}M$
vertices of
$D$ lie in a single component of $G^-\cup H$ tends to $1$  exponentially
rapidly in $n$ (intuitively, that
the sprinkling $H$ joins together the disparate components of
$G^-$).
We will prove this by estimating the complementary probability, which we
will show is in fact much smaller than exponential in $n$.

Suppose that it is not the case that at least $\frac{1}{3}M$
vertices of
$D$ lie in a single component of $G^-\cup H$ (sprinkling fails).
Then there exists $J\subseteq I$ such that the union of $C_j$ over
$j\in J$ has
between $\frac{1}{3}M$ and $\frac{2}{3}M$ vertices of $D$.
 Set $K=I\backslash J$ for convenience and
let $C_J$, $C_K$ denote the union of the $C_i$ over $i\in J$, $i\in K$
respectively.  Then $D=C_J\cup C_K$, and each of $C_J$ and $C_K$
has size between
$\frac{1}{3}M$ and $\frac{2}{3}M$. Fix such $J,K$.  Critically,
there must be no path from $C_J$ to $C_K$ in $H$.

Let $\Delta$ be such that $e^{-\Delta^2/2n}< \frac{1}{6} c_2 \epsilon$,
and set $c_3 = \frac 16 c_2$.
By Lemma~\ref{manypaths},
there are at least $c_3\epsilon 2^nn^{-2\Delta}$ disjoint paths
from $C_J$ to $C_K$ in $\qn$, each of length at most $\Delta$.  Each path is
in $H$ with probability at least $(\epsilon/2n)^{\Delta}$.  Disjointness implies
independence and the probability that $H$ has none of these paths
is at most
\eq \left[1-(\epsilon/2n)^{\Delta}\right]^{c_3\epsilon 2^nn^{-2\Delta}}
\leq \exp\left[-\epsilon^{\Delta}c_3\epsilon 2^nn^{-3\Delta}
2^{-\Delta} \right] . \en

The above quantity bounds the probability that sprinkling fails for
a particular $J,K$.  Since each component $C_i$ is of size at least
$2^{\alpha n/3}$, the number of components $|I|$ is at most
$2^{n(1-\alpha/3)}$.  The number of choices for $J$
(and hence $K=I\backslash J$)
is bounded by $2$ to this number.  Thus the total probability
that sprinkling fails is bounded from above by
\eq
\lbeq{spbd}
    2^{2^{n(1-\alpha /3)}}\exp\left[-\epsilon^{\Delta}c_3\epsilon
    2^nn^{-3\Delta}2^{-\Delta} \right]
    =
    \exp \left[
    (\log 2) 2^{n(1-\alpha /3)}
    -\epsilon^{\Delta}c_3\epsilon
    2^nn^{-3\Delta}2^{-\Delta}
    \right].
\en

Finally, we make some computations to estimate \refeq{spbd}.
Take $\epsilon=e^{-cn^{1/3}}$.  It is at this final stage of the
argument that we use cannot do better than this specific form for $\epsilon$
--- increasing $\epsilon$ only helps.  Then we may take
$\Delta \sim c'n^{2/3}$.  Fix $c_4$ with $1-\alpha/3 < c_4 < 1$.
We select $c$ appropriately small so that $\epsilon^{\Delta}
\geq 2^{-(1-c_4)n}$.  Then $\epsilon^{\Delta}c_3\epsilon
2^nn^{-3\Delta}2^{-\Delta} \geq 2^{(c_4+o(1))n}$, since the factors
$c_3\epsilon n^{-3\Delta}2^{-\Delta}$ are absorbed into the $o(1)$.
Therefore, as required, \refeq{spbd} is exponentially small
(in fact, doubly so).
Thus we have shown that the probability
that sprinkling fails for a particular $J,K$ is much smaller than the
reciprocal of the number $2^{|I|}$ of such $J,K$, and hence a.a.s.\ the
sprinkling succeeds, no $J,K$ exist, and there is a component of size
at least $\frac{1}{3}M$.
\end{proof}

\subsection{The expected cluster size}
\label{sec-chibdsuper}

Finally, we prove the lower bound on $\chi(p)$ stated
in \refeq{chiasysupz}.
Let $\Ccal_1, \Ccal_2, \ldots, \Ccal_{2^n}$ denote the clusters in
$\qn$ arranged in
decreasing order:
    \eq
    |\Ccal_1|= |\Cmax| \geq |\Ccal_2| \geq \ldots ,
    \en
with $\Ccal_i=\varnothing$ if there are fewer than $i$ clusters.
By translation invariance,
    \eqalign
    \lbeq{rewriteclust}
    \chi(p) &= 2^{-n} \sum_{x\in \qn} \expect_p|C(x)|
    =2^{-n} \sum_{x\in \qn}
    \expect_p \Big[\sum_{i=1}^{2^n}|\Ccal_i| \;I[x \in \Ccal_i]\Big]
    \nnb &
    =2^{-n} \expect_p \Big[\sum_{i=1}^{2^n}|\Ccal_i|^2\Big]
    \geq 2^{-n} \expect_p \big[|\Cmax|^2\big].
    \enalign
By \refeq{Cmaxlb},
    \eq
    \chi(p) \geq 2^{-n}[ c_1 \epsilon 2^n]^2
    \Pbold_p(|\Cmax| \geq c_1\epsilon 2^n)
    = [1+\littleo(1)] (c_1 \epsilon)^2 2^n.
    \en
This completes the proof of Theorem~\ref{thm-mr}.

\section*{Acknowledgements}
This work began during a conversation at afternoon tea,
while RvdH, GS and JS were visiting Microsoft Research.
The authors thank Benny Sudakov for bringing the question
of the critical point of the $n$-cube to our attention, and
for telling us about the recent paper \cite{ABS02}.
The work of GS was supported in part by NSERC of Canada.
The work of RvdH was carried out in part at the University
of British Columbia and in part at Delft University of Technology.

\bibliographystyle{plain}

\end{document}

%% file: def99c.tex
\def\UseSection{
        \numberwithin{equation}{section}
    \theoremstyle{plain}
        \newtheorem{theorem}    {Theorem}[section]
        \DefineTheorems 
}

\def\DefineTheorems{
    
    \newtheorem{lemma}      [theorem] {Lemma}
    
    \newtheorem{prop}       [theorem] {Proposition}
    
    \newtheorem{cor}        [theorem] {Corollary}

    \newtheorem{conj}       [theorem] {Conjecture}
    \theoremstyle{definition}
    \newtheorem{defn}       [theorem] {Definition}

    \theoremstyle{definition}

}

\newcommand{\bt}   {\begin{theorem}}
\newcommand{\et}   {\end  {theorem}}
\newcommand{\bl}   {\begin{lemma}}
\newcommand{\el}   {\end  {lemma}}
\newcommand{\bp}   {\begin{prop}}
\newcommand{\ep}   {\end  {prop}}
\newcommand{\bc}   {\begin{cor}}
\newcommand{\ec}   {\end  {cor}}
\newcommand{\bd}   {\begin{defn}}
\newcommand{\ed}   {\end  {defn}}

\newcommand{\ba}   {\begin{array}}
\newcommand{\ea}   {\end  {array}}
\newcommand{\be}   {\begin{enumerate}}
\newcommand{\ee}   {\end  {enumerate}}
\newcommand{\bi}   {\begin{itemize}}
\newcommand{\ei}   {\end  {itemize}}

\def\eq#1\en{\begin{equation}#1\end{equation}}  
\def\eqsplit#1\ensplit{
    \begin{equation}\begin{split}#1\end{split}\end{equation}
    }
\def\eqalign#1\enalign{
    \begin{align}#1\end{align}
    }
\def\eqmul#1\enmul{
    \begin{multline}#1\end{multline}
    }
\newcommand{\eqarrstar} {\begin{eqnarray*}} 
\newcommand{\enarrstar} {\end{eqnarray*}} 
\newcommand{\eqarray}   {\begin{eqnarray}} 
\newcommand{\enarray}   {\end{eqnarray}} 
\newcommand{\nnb}   {\nonumber \\} 

\newcommand{\lbeq}[1]  {\label{e:#1}}
\newcommand{\refeq}[1] {\eqref{e:#1}}    

\makeatletter
\newcommand{\labelcounter}[2]{{%
    \stepcounter{#1}
    \protected@write\@auxout{}%
    {\string\newlabel{#2}{{\csname the#1\endcsname}{\thepage}}}%
    {\ref{#2}}
    }}
\makeatother
\newcommand{\Pbold} {{\mathbb P}}

\newcommand{\Rbold} {{\mathbb R}}

\newcommand{\Zbold} {{\mathbb Z}}

\newcommand{\Ccal}   {\mathcal{C}}

\newcommand{\Zd}    {{ {\Zbold}^d }}

\newcommand{\spose}[1] {{\hbox to 0pt{#1\hss}} }
\newcommand{\ltapprox} {\mathrel{\spose{\lower 3pt\hbox{$\mathchar"218$}}
 \raise 2.0pt\hbox{$\mathchar"13C$}}}
\newcommand{\gtapprox} {\mathrel{\spose{\lower 3pt\hbox{$\mathchar"218$}}
 \raise 2.0pt\hbox{$\mathchar"13E$}}}

